\documentclass[11pt,reqno,oneside]{amsart}
 \usepackage{latexsym}
 \usepackage{amssymb}
 \usepackage{amsfonts}
\usepackage{graphics}
\usepackage{mathcomp}
 
\thanks{University of Li\`ege, Institute of mathematics, Grande Traverse, 12 - B37, B-4000
 Li\`ege, Belgium email :  F.Radoux@ulg.ac.be}
\author[Radoux]{F. Radoux}
\date{\today} 
\title[Explicit formula]{Explicit formula for the natural and projectively equivariant quantization}

\newtheorem{lem}{Lemma}
\newtheorem{thm}[lem]{Theorem}
\newtheorem{prop}[lem]{Proposition}

\theoremstyle{remark}

\theoremstyle{definition}
\newtheorem{defi}{Definition}

\newcommand{\R}{\mathbb{R}}

\newcommand{\N}{\mathbb{N}}

\newcommand{\g}{\mathfrak{g}}

\newcommand{\h}{\mathfrak{h}}

\begin{document}
\begin{abstract}
In \cite{Leconj}, P. Lecomte conjectured the existence of a natural and
projectively equivariant quantization. In \cite{Bord}, M. Bordemann proved
this existence using the framework of Thomas-Whitehead connections. In \cite{moi}, we gave a new proof of the same theorem thanks to
the Cartan connections. After these works, there was no explicit formula for
the quantization. In this paper, we give this formula using the formula in
terms of Cartan connections given in \cite{moi}. This
explicit formula constitutes the generalization to any order of the formulae at second and third
orders soon published by Bouarroudj in \cite{Bou1} and \cite{Bou2}.

\end{abstract}
\maketitle
\section{Introduction}
A quantization can be defined as a linear bijection from the space $\mathcal{S}(M)$ of symmetric contravariant tensor fields
on a manifold $M$ (also called the
space of \emph{Symbols}) to the space $\mathcal{D}_{\frac{1}{2}}(M)$ of
differential operators acting between half-densities.

It is known that
there is no natural quantization procedure. In other words, the spaces of
symbols and of differential operators are not isomorphic as representations of
$\mathrm{Diff}(M)$.

The idea of equivariant quantization, introduced by P. Lecomte and V. Ovsienko
in \cite{LO} is to reduce the group of local diffeomorphisms in the following
way.

They considered the case of the projective group
$PGL(m+1,\R)$ acting locally on the manifold $M=\R^m$ by linear fractional
transformations. They showed that the spaces of symbols and of differential
operators are canonically isomorphic as representations of $PGL(m+1,\R)$ (or
its Lie algebra $sl(m+1,\R)$). In other words, they showed that there exists a
unique \it{projectively equivariant quantization}. \rm{In} \cite{DO}, the authors
generalized this result to the spaces $\mathcal{D}_{\lambda\mu}(\R^{m})$ of
differential operators acting between $\lambda$- and $\mu$-densities and to
their associated graded spaces $\mathcal{S}_{\delta}$. They showed the
existence and uniqueness of a projectively equivariant quantization, provided
the shift value $\delta=\mu-\lambda$ does not belong to a set of critical
values. 

The problem of the $sl(m+1,\R)$-equivariant quantization on $\R^{m}$ has a
counterpart on an arbitrary manifold $M$. In \cite{Leconj}, P. Lecomte conjectured the existence of a quantization 
procedure depending on a torsion-free connection, that would be
 natural (in all arguments) and that would be left invariant by a projective
 change of connection.

After the proof of the existence of such a \it{Natural and equivariant
  quantization} \rm{given} by M. Bordemann in \cite{Bord}, we analysed in \cite{moi} the problem of this existence using
 Cartan connections. We obtained an explicit formula for the quantization map
  in terms of the normal Cartan connection associated to a projective equivalence
 class of torsion free-linear connections. This formula is nothing but the formula for the
 flat case given in \cite{DO} up to replacements of the partial derivatives by
  the invariant differentiation.

The goal of this paper is to obtain an explicit formula on $M$ for the natural
and projectively equivariant quantization. In order to do this, we develop
the operators $\nabla^{\omega^{l}}$ and $Div^{\omega^{l}}$ intervening in the
formula given in \cite{moi} in terms of operators on $M$. This task can be
realized using tools exposed in \cite{Slo1}.

The paper is organized as follows. In the first section, we recall the fundamental notions necessary to
understand the article. In the second part, we calculate the \it{deformation
  tensor}, \rm{the} most important ingredient intervening in the developments
of $\nabla^{\omega^{l}}$ and $Div^{\omega^{l}}$. In the third section, we give
an algorithm that allows to compute these developments thanks to a general
algorithm given in \cite{Slo1}. Finally, in the last part, we calculate the
explicit developments of $\nabla^{\omega^{l}}$ and $Div^{\omega^{l}}$ and we
derive the explicit formula. We show that this formula generalizing the formulae at second and third
orders soon published by Bouarroudj in \cite{Bou1} and \cite{Bou2}.   
\section{Fundamental tools}

For the sake of completeness, we briefly recall in this section the main notions and results
of \cite{moi}. Throughout this note, we denote by $M$ a smooth, Hausdorff and
second countable manifold of dimension $m$.
\subsection{Tensor densities}
The vector bundle of tensor densities $F_{\lambda}(M)\to M$ is a line
bundle associated to the linear frame bundle :
\[ F_{\lambda}(M) = P^1M\times_{\rho}\Delta^{\lambda}(\R^m),\]
where the representation $\rho$ of the group $GL(m,\R)$ on the one-dimensional
vector space $\Delta^{\lambda}(\R^m)$ is given by
\[\rho(A) e = \vert det A\vert^{-\lambda} e,\quad\forall A\in
GL(m,\R),\;\forall e\in \Delta^{\lambda}(\R^m).\]
As usual, we denote by $\mathcal{F}_{\lambda}(M)$ the space of smooth
sections of this bundle. This is the space $C^{\infty}(P^1M,
\Delta^{\lambda}(\R^m))_{GL(m,\R)}$ of functions $f$ such that
\[f(u A) = \rho (A^{-1}) f(u)\quad \forall u \in P^1M,\;\forall A\in
GL(m,\R).\]

\subsection{Differential operators and symbols}
We denote by $\mathcal{D}_{\lambda,\mu}(M)$ the space of
differential operators from $\mathcal{F}_{\lambda}(M)$ to
$\mathcal{F}_{\mu}(M)$. 
The space $\mathcal{D}_{\lambda,\mu}$ is filtered by the order of
differential operators. We denote by $\mathcal{D}^k_{\lambda,\mu}$ the
space of differential operators of order at most $k$. 
The space of \emph{symbols} is then the associated graded space of
$\mathcal{D}_{\lambda,\mu}$.

We denote by $S^l_{\delta}(\R^m)$ the vector space
 $S^l\R^m\otimes\Delta^{\delta}(\R^m)$. There is a natural representation
 $\rho$ of $GL(m,\R)$ on this space (the representation of $GL(m,\R)$ on
 symmetric tensors is the natural one). We then denote by 
$S^l_{\delta}(M)\to M$
 the
vector bundle
\[P^1M\times_{\rho}S^l_{\delta}(\R^m)\to M,\]
and by $\mathcal{S}^l_{\delta}(M)$ the space of smooth sections of
$S^l_{\delta}(M)\to M$, that is, the space $C^{\infty}(P^1M,
S^l_{\delta}(\R^m))_{ GL(m,\R)}$.

 Then if $\delta = \mu
- \lambda$ the \emph{principal symbol operator} $\sigma :
\mathcal{D}^l_{\lambda,\mu}(M)\to \mathcal{S}^l_{\delta}(M)$ commutes with
the action of diffeomorphisms and is a bijection from the quotient space
$\mathcal{D}^l_{\lambda,\mu}(M)/\mathcal{D}^{l-1}_{\lambda,\mu}(M)$ to
$\mathcal{S}^l_{\delta}(M)$.

\subsection{Projective equivalence of connections}
We denote by $\mathcal{C}_M$ the space of torsion-free linear connections
on $M$. Two such connections are \emph{Projectively equivalent} if there
exists a one-form $\alpha$ on $M$ such that their associated covariant
derivatives $\nabla$ and $\nabla'$ fulfill the relation
\[\nabla'_XY = \nabla_XY + \alpha(X) Y + \alpha(Y)X.\]
 
\subsection{Problem setting}
A \emph{quantization} on $M$ is a linear bijection $Q_M$ from
 the space of symbols
$\mathcal{S}_{\delta}(M)$ to the space of differential operators
$\mathcal{D}_{\lambda,\mu}(M)$ such that
\[\sigma (Q_M(S)) = S,\quad\forall S\in \mathcal{S}^k_{\delta}(M),\; \forall k
\in\mathbb{N}. \]
A \emph{natural quantization} is a quantization which
depends on a torsion-free connection and commutes with
 the action of diffeomorphisms (see \cite{moi} for a more precise definition).

A quantization $Q_M$ is \emph{projectively equivariant} if one has
$Q_M(\nabla) = Q_M(\nabla')$ whenever $\nabla$ and $\nabla'$ are projectively
equivalent torsion-free linear connections on $M$.

\subsection{Projective structures and Cartan projective connections}
These tools were presented in detail in \cite[Section 3]{moi}. We give here the
 most important ones for this paper to be self-contained.

 We consider the group $G=PGL(m+1,\R)$.
 We denote by $H$ 
the subgroup  
\begin{equation}\label{hproj}H=\{\left(\begin{array}{cc}A & 0\\\xi & a
\end{array}\right) : A\in GL(m,\R),\xi\in \R^{m*}, a\not = 0\}/\R_0\mbox{Id}.
\end{equation}
The group $H$ is the semi-direct product $G_0 \rtimes G_1$, where
$G_0$ is isomorphic to $GL(m,\R)$ and $G_1$ is isomorphic to $\R^{m*}$. The
Lie algebra associated to $H$ is $\g_0\oplus\g_1$. 

It is well-known that $H$ can be seen as a subgroup of the group of 2-jets 
$G_m^2$.\\
A \emph{Projective
  structure on $M$} is then a reduction of the second order frame bundle $P^2M$
to the group $H$. \\
The following result (\cite[p. 147]{Koba}) is 
the starting point of our method :
\begin{prop}[Kobayashi-Nagano]
There is a natural one to one correspondence between the projective
 equivalence classes of torsion-free linear connections on $M$
and the projective structures on $M$.
\end{prop}
We now recall the definition of a projective Cartan connection :
\begin{defi}
Let $P\to M$ be a principal $H$-bundle. A projective Cartan connection on $P$
is a $sl(m+1,\R)$- valued 1-form $\omega$ such that
\begin{itemize}
\item There holds $R_a^*\omega=Ad(a^{-1})\omega,\quad\forall a \in H$,
\item One has $\omega(k^*)=k\quad\forall k\in \h=\g_0\oplus \g_1$,
\item For all $u\in P$, $\omega_u :T_uP\to sl(m+1,\R)$ is a linear bijection.
\end{itemize}
\end{defi}
In general, if $\omega$ is a Cartan connection defined on a $H$-principal bundle $P$, then its
curvature $\Omega$ is defined as usual by
 \begin{equation}\label{curv} 
\Omega = d\omega+\frac{1}{2}[\omega,\omega].
\end{equation}
We can define from $\Omega$ a function $\kappa\in
C^{\infty}(P,\g_{-1}^{*}\otimes\g_{-1}^{*}\otimes\g)$ by :
$$\kappa(u)(X,Y):=\Omega(u)(\omega^{-1}(X),\omega^{-1}(Y)).$$
The \emph{Normal} Cartan connection has the following property (see
\cite[p. 136]{Koba}):
$$\sum_{i}\kappa_{jil}^{i}=0\quad\forall j,\forall l.$$

Now, the following result (\cite[p. 135]{Koba}) gives the relationship between  projective
structures and Cartan connections :
\begin{prop}
 A unique normal
 Cartan projective connection is
 associated to every projective structure $P$. This association is natural.
\end{prop}
The connection associated to a projective structure $P$ is called the normal
projective connection of the projective structure.
\subsection{Lift of equivariant functions}\label{Lift}

If $(V,\rho)$ is a representation of $\mathit{GL}(m,\R)$, then we can define
from it a representation $(V,\rho')$ of $H$ (see \cite{moi} section 3). If $P$ is  a projective structure on $M$, the natural projection $P^2M\to
P^1M$ induces a projection $p :P\to P^1M$ and we have a well-known result:
\begin{prop} If $(V,\rho)$ is a representation of $GL(m,\R)$, then the map
\[p^* : C^{\infty}(P^1M,V)\to C^{\infty}(P,V) : f \mapsto f\circ p\]
defines a bijection from $C^{\infty}(P^1M,V)_{\mathrm{GL}(m,\R)}$ to
$C^{\infty}(P,V)_{H}$.
\end{prop}

Subsequently, we will use the representation $\rho'_*$ of the Lie algebra of
$H$ on $V$. If we recall that this algebra is isomorphic to
 $gl(m,\R)\oplus \R^{m*}$
then we have 
\begin{equation}\label{rho}\rho'_* (A, \xi) = \rho_*(A),\quad\forall A\in
  gl(m,\R), \xi\in \R^{m*}.\end{equation}
Recall that if $f\in C^{\infty}(P,V)_H$ then one has
\begin{equation}\label{Invalg}
L_{h^*}f(u) + \rho'_*(h)f(u)=0,\quad\forall h\in gl(m,\R)\oplus\R^{m*}\subset sl(m+1,\R),
\forall u\in P.
\end{equation}
\subsection{The first explicit formula}
First, we give the definitions of operators used subsequently :
\begin{defi}
Let $(V, \rho)$ be a representation of $H$. If $f\in C^{\infty}(P,V)$, then $\nabla^{{\omega}^k} f
\in C^{\infty}(P,\otimes^k\R^{m*}\otimes V)$ is defined by
\[\nabla^{{\omega}^k} f(u)(X_1,\ldots,X_k) = L_{\omega^{-1}(X_{k})}\circ\ldots\circ L_{\omega^{-1}(X_{1})}f(u).\]
\end{defi}
If we symmetrize this operation, we obtain the
\begin{defi}
If $f\in C^{\infty}(P,V)$, then $\nabla_s^{{\omega}^k} f
\in C^{\infty}(P,S^k\R^{m*}\otimes V)$ is defined by :
\[\nabla_s^{{\omega}^k} f(u)(X_1,\ldots,X_k) = \frac{1}{k!}\sum_{\nu}\nabla^{{\omega}^k} f(u)(X_{\nu_1},\ldots,X_{\nu_k}).\]
\end{defi}
If $(e_1,\ldots,e_m)$ is the canonical basis of $\R^m$ and if
$(\epsilon^1,\ldots,\epsilon^m)$ is the dual basis corresponding in $\R^{m*}$,
the \emph{divergence operator} is defined then by :
\[Div^{\omega} : C^{\infty}(P,S^k_{\delta}(\R^m))
\to C^{\infty}(P,S^{k-1}_{\delta}(\R^m)) :
S\mapsto \sum_{j=1}^m \nabla^{\omega}_{e_j}S(\epsilon^{j}).\]

If $\gamma\in C^{\infty}(P^{1}M, \Delta^{\lambda}(\R^{m})\otimes
S^{l}\R^{m*})$, one defines the symmetrized covariant derivative of $\gamma$, $\nabla_s\gamma\in C^{\infty}(P^{1}M, \Delta^{\lambda}(\R^{m})\otimes
S^{l+1}\R^{m*})$, by :
$$(\nabla_s\gamma)(X_1,\ldots,X_{l+1})=\frac
{1}{(l+1)!}\sum_{\nu}(\nabla_{X_{\nu(1)}}\gamma)(X_{\nu(2)},\ldots,X_{\nu(l+1)}).$$
Recall now the definition of the numbers $\gamma_{2k-l}$ :
\[\gamma_{2k-l} = \frac{m+2 k - l -(m+1)\delta}{m+1}.\]
A value of $\delta$ is \it{critical} \rm{if there are} $k,l\in \N$ such that
$1\leq l\leq k$ and $\gamma_{2k-l}=0$.

Finally, we can recall the formula giving the natural and projectively
equivariant quantization in terms of the normal Cartan connection (see \cite{moi}, theorem 11) :
\begin{thm} If $\delta$ is not critical, then the
 collection of maps \\
$Q_M : \mathcal{C}_M\times \mathcal{S}_{\delta}(M)\to
\mathcal{D}_{\lambda,\mu}(M)$ defined by
\begin{equation}\label{formula}Q_M(\nabla, S)(f) = p^{*^{-1}}(\sum_{l=0}^k C_{k,l} \langle Div^{\omega^l}
p^*S,\nabla_s^{\omega^{k-l}}p^*f\rangle),\forall S\in \mathcal{S}^k_{\delta}(M)\end{equation}
defines a projectively invariant natural quantization if 
\[C_{k,l} =\frac{(\lambda + \frac{k-1}{m+1})\cdots (\lambda +
  \frac{k-l}{m+1})}{\gamma_{2k-1}\cdots
  \gamma_{2k-l}}\left(\begin{array}{c}k\\l\end{array}\right),\forall l\geq 1,\quad C_{k,0}=1.\]
\end{thm}
\section{The deformation tensor}
An Ehresmann connection $\gamma$ on $P^{1}M$ belonging to a projective
structure $P$ induces a $GL(m,\R)$-equivariant section $\sigma$ of $P\to
P^{1}M$ (see \cite{Koba} page 147). This correspondence
establishes a bijection between the set of the connections belonging to the
projective structure $P$ and the set of the $GL(m,\R)$-equivariant sections
of $P\to P^{1}M$. 

\vspace{0.3cm}\rm{If} $\sigma$ is the section corresponding to a connection
$\gamma$, one can define an application $\tau:P\to\g_1$ in
the following way :
$$u=\sigma(p(u)).\exp(\tau(u)).$$

\vspace{0.2cm}If $\gamma$ is a connection on $P^{1}M$ corresponding to a
section $\sigma$ and if $\omega$ is the normal Cartan connection corresponding to the projective
class of $\gamma$, one has the following result (see \cite{Slo1} page 43) :
\begin{prop}
There is a
unique Cartan connection $\tilde{\gamma}=\omega_{-1}\oplus\omega_{0}\oplus\tilde{\gamma}_{1}$
such as $\tilde{\gamma}_{1}|(T\sigma(TP^{1}M))=0.$
\end{prop}
This Cartan connection is called the \it{Cartan connection induced
  by $\gamma$}.

\vspace{0.2cm}\rm{The normal} Cartan connection $\omega$ and the Cartan connection $\tilde{\gamma}$ induced by
$\gamma$ differ only by their components in $\g_1$. Moreover, as the
difference $\omega-\tilde{\gamma}$ vanishes on vertical vector fields, there
is a function $\Gamma\in C^{\infty}(P,\g_{-1}^{*}\otimes\g_1)$ such as
$$\omega=\tilde{\gamma}-\Gamma\circ\omega_{-1}.$$
This function is $H$-equivariant and represents then a tensor of type
$\left(\begin{array}{c}0\\2\end{array}\right)$ on $M$ ;
it is called the \it{deformation tensor} \rm{(see \cite{Slo1} page 45)}. This
function has the following property (see \cite{Slo1} lemma 3.10) :
\begin{equation}\label{def}
(\tilde{\kappa}_{0}-\kappa_{0})(u)(X,Y)=[X,\Gamma(u).Y]+[\Gamma(u).X,Y]
\end{equation}

\vspace{0.2cm}if $u\in P$, $X,Y\in\g_{-1}$ and if $\tilde{\kappa}_{0}$ and
$\kappa_{0}$ are the functions induced respectively by the curvatures of
$\tilde{\gamma}$ and of $\omega$.

One can compute the deformation tensor in the projective case exactly in the
same way as it is calculated in the conformal case at page 63 of \cite{Slo1}.
First we fix a basis $e_{i}$ in $\g_{-1}$, $e_{j}^{i}$ in $\g_0$, $\epsilon^{i}$ in
$\g_1$. We have then
$$\Gamma(u)(e_i)=\sum_{j}\Gamma(u)_{ji}\epsilon^{j},$$
$$\kappa_{0}(u)(e_i,e_j)=\sum_{k,l}\kappa_{0}(u)_{lij}^{k} e_{k}^{l}$$
and
$$\tilde{\kappa}_{0}(u)(e_i,e_j)=\sum_{k,l}\tilde{\kappa_{0}}(u)_{lij}^{k}
e_{k}^{l}.$$

One obtains then using the equality (\ref{def}) the following relations :
\begin{equation}\label{defo1}
(\kappa_0-\tilde{\kappa_0})_{klj}^{l}=\Gamma_{jk}-m\Gamma_{kj}; 
\end{equation}
\begin{equation}\label{defo2}
(\kappa_0-\tilde{\kappa_0})_{kij}^{k}=(m+1)(\Gamma_{ji}-\Gamma_{ij}).
\end{equation}

\vspace{0.2cm}On one hand, the functions $(\kappa_{0})_{klj}^{l}$ and
$(\kappa_{0})_{kij}^{k}$ vanish by normality of $\omega$. On the other hand,
the functions $(\tilde{\kappa_{0}})_{lij}^{k}$ are the components of the
equivariant function on $P$ that represents the curvature tensor corresponding to the connection $\gamma$. A straightforward computation allows then to obtain the
expression of the deformation tensor from the relations (\ref{defo1}) and (\ref{defo2}) :
\begin{equation}\label{Gamma}
\Gamma_{jk}=\frac{\mathrm{Ric}_{kj}}{1-m}+\frac
      {m\;\mathrm{trR}_{jk}}{(m+1)(m-1)},
\end{equation}
where $\mathrm{Ric}$ and $\mathrm{trR}$ represent the equivariant functions on
$P$ corresponding respectively to the Ricci tensor and to the trace of the curvature. 
\section{Developments of $\nabla^{\omega^{l}}$ and $Div^{\omega^{l}}$}
In order to obtain an explicit formula for the quantization, we need to know the
developments of the operators $\nabla^{\omega^{l}}$ and $Div^{\omega^{l}}$ in terms of
operators on $M$. We first recall the developments of \cite{Slo1}.

Let $\gamma$ be an Ehresmann connection on $P^{1}M$ corresponding to a
covariant derivative $\nabla$ and belonging to a
projective structure $P$. We denote by $\omega$ the normal Cartan connection
on $P$.

Let $(V,\rho)$ be a representation of $GL(m,\R)$ inducing a representation
$(V,\rho_{*})$ of $gl(m,\R)$. If we denote by $\rho_{*}^{(l)}$ the canonical
representation on $\otimes^{l}\g_{-1}^{*}\otimes V$ and if $s\in C^{\infty}(P^{1}M,V)_{GL(m,\R)}$, then
$F^{l}s:=\nabla^{\omega^{l}}(p^{*}s)-p^{*}(\nabla^{l}s)$ is given by the
following induction :

\vspace{0.2cm}$F^{0}s(u)=0$

\vspace{0.2cm}$F^{l}s(u)(X_1,\ldots,X_l)=\rho_{*}^{(l-1)}([X_{l},\tau(u)])(F^{l-1}s(u))(X_1,\ldots,X_{l-1})$

\vspace{0.2cm}\hspace{2cm}$+S_{\tau}(F^{l-1}s(u))(X_1,\ldots,X_{l-1})$

\vspace{0.2cm}\hspace{2cm}$+S_{\nabla}(F^{l-1}s(u))(X_1,\ldots,X_{l-1})$

\vspace{0.2cm}\hspace{2cm}$+S_{\Gamma}(F^{l-1}s(u))(X_1,\ldots,X_{l-1})$

\vspace{0.2cm}\hspace{2cm}$+\rho_{*}^{(l-1)}([X_{l},\tau(u)])(p^{*}(\nabla^{l-1}s)(u))(X_1,\ldots,X_{l-1}).$

\vspace{0.5cm}This expression expands into a sum of terms of the form
$$a\rho_{*}^{(t_1)}(\beta_1)\ldots\rho_{*}^{(t_i)}(\beta_{i})p^{*}\nabla^{j}s$$
where $a$ is a scalar coefficient, the $\beta_{l}$ are iterated brackets
involving some arguments $X_{l}$, the iterated invariant differentials
$\nabla^{r}\Gamma$ evaluated on some arguments $X_{l}$, and $\tau$. Exactly
the first $t_j$ arguments $X_1,\ldots,X_{t_j}$ are evaluated after the action
of $\rho_{*}^{(t_j)}(\beta_j)$, the other ones appearing on the right are
evaluated before. The individual transformations in $S_{\tau}$, $S_{\nabla}$
and $S_{\Gamma}$ act as follows : 
\begin{enumerate}
\item The action of $S_{\tau}$ replaces each summand
  $a\rho_{*}^{(t_1)}(\beta_1)\ldots\rho_{*}^{(t_i)}(\beta_{i})p^{*}\nabla^{j}$
  by a sum with just one term for each occurrence of $\tau$ where this $\tau$
  is replaced by $[\tau,[\tau,X_{l}]]$ and the coefficient $a$ is multiplied
  by $-\frac {1}{2}$.
\item The transformation $S_{\nabla}$ replaces each summand in $F^{l-1}$ by a
  sum with just one term for each occurrence of $\Gamma$ and its
  differentials, where these arguments are replaced by their covariant
  derivatives $\nabla_{X_l}$, and with one additional term where $\nabla^{j}s$
  is replaced by $\nabla_{X_l}(\nabla^{j}s)$. 
\item The transformation $S_{\Gamma}$ replaces each summand by a sum with just
  one term for each occurrence of $\tau$ where this $\tau$ is replaced by $\Gamma(u).X_{l}$.
\end{enumerate} 
In fact, this algorithm becomes easily linear in the following way :
\begin{prop}
The development of $\nabla^{\omega^{l}}(p^{*}s)(X_1,\ldots,X_{l})$ is obtained
as follows :

\vspace{0.2cm}$\nabla^{\omega^{l}}(p^{*}s)(X_1,\ldots,X_l)=\rho_{*}^{(l-1)}([X_{l},\tau])(\nabla^{\omega^{l-1}}(p^{*}s))(X_1,\ldots,X_{l-1})$

\vspace{0.2cm}\hspace{2cm}$+S_{\tau}(\nabla^{\omega^{l-1}}(p^{*}s))(X_1,\ldots,X_{l-1})$

\vspace{0.2cm}\hspace{2cm}$+S_{\nabla}(\nabla^{\omega^{l-1}}(p^{*}s))(X_1,\ldots,X_{l-1})$

\vspace{0.2cm}\hspace{2cm}$+S_{\Gamma}(\nabla^{\omega^{l-1}}(p^{*}s))(X_1,\ldots,X_{l-1}).$
\end{prop}
\begin{prop}
If $f\in C^{\infty}(P^{1}M,\Delta^{\lambda}(\R^{m}))_{GL(m,\R)}$, then
$\nabla^{\omega^{l}}(p^{*}f)(X,\ldots,X)$ is a linear combination of terms of the form
$$(\otimes^{n_{-1}}\tau\otimes p^{*}(\otimes^{n_{l-2}}\nabla^{l-2}\Gamma\otimes\ldots\otimes\otimes^{n_0}\Gamma\otimes\nabla^{q}f))(X,\ldots,X).$$
If we denote by $T(n_{-1},\ldots,n_{l-2},q)$ a such term,
$T(n_{-1},\ldots,n_{l-2},q)$ gives rise in the development of
$\nabla^{\omega^{l+1}}(p^{*}f)(X,\ldots,X)$ to
$$(-\lambda(m+1)-2l+n_{-1})T(n_{-1}+1,\ldots,n_{l-2},q)+T(n_{-1},\ldots,n_{l-2},q+1)$$
$$+\sum_{j=-1}^{l-2}n_{j}T(n_{-1},\ldots,n_{j}-1,n_{j+1}+1,\ldots,n_{l-2},q).$$
\end{prop}
\begin{proof}
One sees indeed easily that the application of the first part of the algorithm
gives 
$$(-\lambda(m+1)-2l)T(n_{-1}+1,\ldots,n_{l-2},q).$$
The second part gives
$$n_{-1}T(n_{-1}+1,\ldots,n_{l-2},q).$$
The third part contributes to 
$$T(n_{-1},\ldots,n_{l-2},q+1)+\sum_{j=0}^{l-2}n_{j}T(n_{-1},\ldots,n_{j}-1,n_{j+1}+1,\ldots,n_{l-2},q).$$
The fourth
gives 
$$n_{-1}T(n_{-1}-1,n_{0}+1,\ldots,n_{l-2},q).$$
\end{proof} 
One deduces easily from this result the following corollary :
\begin{prop}\label{devnabla}
If $f\in C^{\infty}(P^{1}M,\Delta^{\lambda}(\R^{m}))_{GL(m,\R)}$, $\nabla_s^{\omega^{l}}(p^{*}f)$ is a linear combination
of terms of the form 
$$\tau^{n_{-1}}\vee p^{*}((\nabla_s^{l-2}r)^{n_{l-2}}\vee\ldots\vee
r^{n_0}\vee\nabla_s^{q}f),$$
where $r$ denotes the symmetric part of the Ricci tensor divided by $1-m$. If
we denote by $T(n_{-1},\ldots,n_{l-2},q)$ a such term,
$T(n_{-1},\ldots,n_{l-2},q)$ gives rise in the development of
$\nabla_s^{\omega^{l+1}}(p^{*}f)$ to
$$(-\lambda(m+1)-2l+n_{-1})T(n_{-1}+1,\ldots,n_{l-2},q)+T(n_{-1},\ldots,n_{l-2},q+1)$$
$$+\sum_{j=-1}^{l-2}n_{j}T(n_{-1},\ldots,n_{j}-1,n_{j+1}+1,\ldots,n_{l-2},q).$$ \end{prop}
\begin{proof}
First we remark that the symmetric part of $\Gamma$ is reduced to $r$ by
antisymmetry of the tensor $\mathrm{tr R}$. It suffices then to remark that if \\
$\nabla^{\omega^{l}}(p^{*}f)(X,\ldots,X)$ is equal to a linear combination of terms of the form
$$(\otimes^{n_{-1}}\tau\otimes
p^{*}(\otimes^{n_{l-2}}\nabla^{l-2}\Gamma\otimes\ldots\otimes\otimes^{n_0}\Gamma\otimes\nabla^{q}f))(X,\ldots,X),$$
then $\nabla_s^{\omega^{l}}(p^{*}f)$ is equal to the
corresponding linear combination of the terms of the form 
$$\tau^{n_{-1}}\vee p^{*}((\nabla_s^{l-2}r)^{n_{l-2}}\vee\ldots\vee
r^{n_0}\vee\nabla_s^{q}f).$$
Indeed, the two last tensors are then equal because they are both symmetric
and that they are equal when they are evaluated on $X^{l}$.
\end{proof}
Remark that the action of the algorithm on the generic term of the development
of $\nabla_s^{\omega^{l}}(p^{*}f)$ can be summarized. Indeed, this action
gives first
$$(-\lambda(m+1)-2l+n_{-1})T(n_{-1}+1,\ldots,n_{l-2},q).$$
It gives next
$$n_{-1}T(n_{-1}-1,n_{0}+1,\ldots,n_{l-2},q).$$
Finally, it makes act the covariant derivative $\nabla_s$ on
$$(\nabla_s^{l-2}r)^{n_{l-2}}\vee\ldots\vee
r^{n_0}\vee\nabla_s^{q}f.$$
\begin{prop}\label{Div}
If $S\in C^{\infty}(P^{1}M,\Delta^{\delta}\R^{m}\otimes
S^{k}\R^{m})_{GL(m,\R)}$, then $Div^{\omega^l}(p^{*}S)$ is a
linear combination of terms of the form
$$\langle\tau^{n_{-1}}\vee p^{*}((\nabla_s^{k-2}r)^{n_{k-2}}\vee\ldots\vee
r^{n_0}), p^{*}(Div^{q}S)\rangle.$$
If we denote by $T(n_{-1},\ldots,n_{l-2},q)$ a such term,
$T(n_{-1},\ldots,n_{l-2},q)$ gives rise in the development of
$Div^{\omega^{l+1}}(p^{*}S)$ to
$$(\gamma_{2(k-l)-1}(m+1)+n_{-1})T(n_{-1}+1,\ldots,n_{l-2},q)+T(n_{-1},\ldots,n_{l-2},q+1)$$
$$+\sum_{j=-1}^{l-2}n_{j}T(n_{-1},\ldots,n_{j}-1,n_{j+1}+1,\ldots,n_{l-2},q).$$
\end{prop}
\begin{proof}
We have to compute
$$(\nabla^{\omega^{l+1}}(p^{*}S)(e_{i_1},\ldots,e_{i_{l+1}}))(\epsilon^{i_1},\ldots,\epsilon^{i_{l+1}}).$$
As the first part of the development of
$$\nabla^{\omega^{l+1}}(p^{*}S)(e_{i_1},\ldots,e_{i_{l+1}})$$
according to the algorithm is
$$(\rho_{*}^{(l)}([e_{i_{l+1}},\tau(u)])\nabla^{\omega^{l}}(p^{*}S)(u))(e_{i_1},\ldots,e_{i_l}),$$
we must first calculate
$$[(\rho_{*}^{(l)}([e_{i_{l+1}},\tau(u)])\nabla^{\omega^{l}}(p^{*}S)(u))(e_{i_1},\ldots,e_{i_l})](\epsilon^{i_1},\ldots,\epsilon^{i_{l+1}}).$$
This latter expression is equal to
$$[\rho_{*}([e_{i_{l+1}},\tau(u)])(\nabla^{\omega^{l}}(p^{*}S)(u)(e_{i_1},\ldots,e_{i_l}))](\epsilon^{i_1},\ldots,\epsilon^{i_{l+1}})$$ 
$$-\sum_{j=1}^{l}(\nabla^{\omega^{l}}(p^{*}S)(u)(e_{i_1},\ldots,[e_{i_{l+1}},\tau(u)]e_{i_j},\ldots,e_{i_l}))(\epsilon^{i_1},\ldots,\epsilon^{i_{l+1}}),$$
i.e. to
$$[\rho'_{*}([e_{i_{l+1}},\tau(u)])(\nabla^{\omega^{l}}(p^{*}S)(u)(e_{i_1},\ldots,e_{i_l})(\epsilon^{i_1},\ldots,\epsilon^{i_{l}}))](\epsilon^{i_{l+1}})$$
$$+\sum_{j=1}^{l}(\nabla^{\omega^{l}}(p^{*}S)(u)(e_{i_1},\ldots,e_{i_l}))(\epsilon^{i_1},\ldots,\epsilon^{i_j}[e_{i_{l+1}},\tau(u)],\ldots,\epsilon^{i_{l+1}})$$
$$-\sum_{j=1}^{l}(\nabla^{\omega^{l}}(p^{*}S)(u)(e_{i_1},\ldots,[e_{i_{l+1}},\tau(u)]e_{i_j},\ldots,e_{i_l}))(\epsilon^{i_1},\ldots,\epsilon^{i_{l+1}}),$$
 if $\rho'$ denotes the action of $GL(m,\R)$ on symbols of degree $k-l$. The
 second and third lines of the previous expression give respectively $2l$ and $-2l$
 terms in which $n_{-1}$ is replaced by $n_{-1}+1$. Their contributions
 vanish. One sees easily that the first line gives
$$\gamma_{2(k-l)-1}(m+1)T(n_{-1}+1,\ldots,n_{l-2},q).$$
One can see that the substitutions intervening in the third last parts of the
algorithm ``commute'' with the valuations in
$\epsilon^{i_1},\ldots,\epsilon^{i_{l}}$ thanks to the general form of
$\nabla^{\omega^{l}}(p^{*}S)(X_{1},\ldots,X_{l})$. Indeed, one can show easily
that $\nabla^{\omega^{l}}(p^{*}S)(X_{1},\ldots,X_{l})$ is a linear combination of terms constructed in
the following way. One applies first $p^{*}(\nabla^{q}S)$ on some $X_{i}$ and
one contracts the result several times with $\tau$. One contracts then the
obtained symbol with tensors of degree 1 obtained contracting some
$p^{*}(\nabla^{t}\Gamma)$ with $t+1$ arguments $X_{i}$. One multiplies symmetrically the result by others $X_i$. Finally, one
multiplies the result by numbers obtained applying $\tau$ on some $X_{i}$ and
some $p^{*}(\nabla^{t}\Gamma)$ on $t+2$ arguments $X_{i}$. 

One sees then that the second part of the algorithm gives $n_{-1}$ terms where $n_{-1}$ becomes
$n_{-1}+1$. One sees too that the third part contributes to
$$T(n_{-1},\ldots,n_{l-2},q+1)+\sum_{j=0}^{l-2}n_{j}T(n_{-1},\ldots,n_{j}-1,n_{j+1}+1,\ldots,n_{l-2},q).$$
The fourth gives 
$$n_{-1}T(n_{-1}-1,n_{0}+1,\ldots,n_{l-2},q).$$
\end{proof}
Remark that the action of the algorithm on the generic term of the development
of $Div^{\omega^{l}}(p^{*}S)$ can be summarized. Indeed, this action gives
first 
$$(\gamma_{2(k-l)-1}(m+1)+n_{-1})T(n_{-1}+1,\ldots,n_{l-2},q).$$
It gives next 
$$n_{-1}T(n_{-1}-1,n_{0}+1,\ldots,n_{l-2},q).$$
 Finally, it makes act the divergence $Div$ on 
$$\langle(\nabla_s^{k-2}r)^{n_{k-2}}\vee\ldots\vee
r^{n_0}, Div^{q}S\rangle.$$
\section{The main result}
Because of the previous propositions, the quantization can be written as a
linear combination of terms of the form
$$\langle\langle\tau^{n_{-1}}\vee p^{*}((\nabla_s^{k-2}r)^{n_{k-2}}\vee\ldots\vee
r^{n_0}), p^{*}(Div^{q}S)\rangle, p^{*}(\nabla_s^{l}f)\rangle.$$
In this expression, it suffices to consider the terms for which
$n_{-1}=0$. Indeed, suppose that the expression
\begin{equation}\label{expr}
\sum_{j=0}^{k}\langle a_{j}, \tau^{j}\rangle
\end{equation}
in which the functions $a_{j}$ are $H$-equivariant is $H$-equivariant. First
note that $L_{h^{*}}\tau=h$ for all
$h\in\g_{1}$ (see
\cite{Slo1} page 48). The fact that the application of $L_{h^{*}}$ to (\ref{expr})
gives $0$ for all $h\in\g_{1}$ implies that 
$$\sum_{j=1}^{k}\langle ja_{j}, \tau^{j-1}\rangle$$
is equal to zero, hence $H$-equivariant. Repeating the process, one finds
finally that
$a_k=0$. One deduces then progressively that the functions $a_{j}$ are equal to zero for $j$ equal to $1,\ldots,k$.  

The following results give the explicit developments of
$\nabla_s^{\omega^{l}}(p^{*}f)$ and of $Div^{\omega^{l}}(p^{*}S)$ :
\begin{prop}
The term of degree $t$ in $\tau$ in the development of
$\nabla_s^{\omega^{l}}(p^{*}f)$ is equal to
$$\left(\begin{array}{c}l\\t\end{array}\right)\prod_{j=1}^{t}(-\lambda(m+1)-l+j)p^{*}(\pi_{l-t}(\sum_{j=0}^{l-t}(\nabla_s+T_{1})^{j})f),$$
if $\pi_{l-t}$ denotes the projection on the operators of degree $l-t$ and if
the restriction of $T_{1}$ to tensors of type
$\left(\begin{array}{c}0\\j\end{array}\right)$ with values
in the $\lambda$-densities is equal to 
$$(-\lambda(m+1)-j)(j+1)$$
times the symmetric product by $r$ (the degree of $\nabla_{s}$ is $1$ whereas
the degree of $T_{1}$ is $2$). We set that $\prod_{j=1}^{t}(-\lambda(m+1)-l+j)$ is equal to $1$ if $t=0$.  
\end{prop} 
\begin{proof}
In order to simplify the notations, denote by $\beta$ the number $-\lambda(m+1)$. The
formula is true if $l$ and $t$ are equal to $0$. Suppose that the formula is
satisfied for all $t$ until the order $l-1$. If
$l-t\geq 2$ and if $t\geq 2$, then the term of degree $t$ in $\tau$ at the
order $l$ is equal using the induction procedure to :
$$(t+1)\left(\begin{array}{c}l-1\\t+1\end{array}\right)\prod_{j=1}^{t+1}(\beta-l+1+j)p^{*}(r\vee\pi_{l-t-2}(\sum_{j=0}^{l-t-2}(\nabla_s+T_{1})^{j})f)$$
$$+\left(\begin{array}{c}l-1\\t\end{array}\right)\prod_{j=1}^{t}(\beta-l+1+j)p^{*}(\nabla_s(\pi_{l-t-1}(\sum_{j=0}^{l-t-1}(\nabla_s+T_{1})^{j}))f)$$
$$+\left(\begin{array}{c}l-1\\t-1\end{array}\right)(\prod_{j=1}^{t-1}(\beta-l+1+j))(\beta-2l+t+1)$$
$$p^{*}(\pi_{l-t}(\sum_{j=0}^{l-t}(\nabla_s+T_{1})^{j})f).$$ 
Note that
$$(\beta-l+t+2)(l-t-1)r\vee\pi_{l-t-2}(\sum_{j=0}^{l-t-2}(\nabla_s+T_{1})^{j})$$
is equal to
$$\pi_{l-t}(T_{1}(\sum_{j=0}^{l-t-2}(\nabla_s+T_{1})^{j})).$$
The sum of the three terms above is then equal to a multiple of 
$$p^{*}(\pi_{l-t}(\sum_{j=0}^{l-t}(\nabla_s+T_{1})^{j})f),$$
this multiple being equal to
$$\prod_{j=2}^{t}(\beta-l+j)(\left(\begin{array}{c}l-1\\t\end{array}\right)(\beta-l+t+1)+\left(\begin{array}{c}l-1\\t-1\end{array}\right)(\beta-2l+t+1)),$$
i.e. to
$$\prod_{j=2}^{t}(\beta-l+j)$$
$$((\beta-l+1)(\left(\begin{array}{c}l-1\\t\end{array}\right)+\left(\begin{array}{c}l-1\\t-1\end{array}\right))+t\left(\begin{array}{c}l-1\\t\end{array}\right)+(t-l)\left(\begin{array}{c}l-1\\t-1\end{array}\right)).$$

\vspace{0.2cm}We conclude using the formula of the Pascal's triangle.

\vspace{0.2cm}We deal with the cases $\;l-t\geq 2\;\&\;t<2,\quad
l-t<2\;\&\;t\geq 2$ and $l-t<2\;\&\;t<2$ in a same way.
\end{proof}

\begin{prop}
The term of degree $t$ in $\tau$ in the development of
$Div^{\omega^{l}}(p^{*}S)$ is equal to 
$$\left(\begin{array}{c}l\\t\end{array}\right)\prod_{j=1}^{t}(\gamma_{2k-1}(m+1)-l+j)p^{*}(\pi_{l-t}(\sum_{j=0}^{l-t}(Div+T_{2})^{j})S),$$
if $\pi_{l-t}$ denotes the projection on the operators of degree $l-t$ and if
the restriction of $T_{2}$ to symbols of degree $j$ is equal to
$$((m+1)\gamma_{2k-1}-k+j)(k-j+1)$$
times the inner product by $r$ (the degree of $Div$ is $1$ whereas the degree
of $T_2$ is $2$). We set that the product \\
$\prod_{j=1}^{t}(\gamma_{2k-1}(m+1)-l+j)$ is equal to $1$ if $t=0$. 
\end{prop} 
\begin{proof}
The proof is completely similar to the one of the previous proposition.
\end{proof}
We can now write the explicit formula giving the natural projectively
equivariant quantization of \cite{moi} :
\begin{thm}
The quantization $Q_{M}$ of \cite{moi} is given by the following \\
formula :
$$Q_{M}(\nabla,S)(f)=\sum_{l=0}^{k}C_{k,l}\langle\pi_{l}(\sum_{j=0}^{l}(Div+T_{2})^{j})S,
\pi_{k-l}(\sum_{j=0}^{k-l}(\nabla_s+T_{1})^{j})f\rangle.$$
\end{thm}
 
\vspace{0.2cm}One can easily derive from this formula the formula at the third order given
by Bouarroudj in \cite{Bou2}. Indeed, if we denote by $D$, $T$, $\partial
T$ the operators $\nabla_{s}$, $r\vee$ and $(\nabla_{s}r)\vee$ (resp. $Div$,
$i(r)$ and $i(\nabla_{s}r)$) and if we denote by $\beta$ the number
$-\lambda(m+1)$ (resp. $\gamma_{5}(m+1)$), one obtains :
$$\pi_{1}(\sum_{j=0}^{1}(D+T)^{j})=D,\quad
\pi_{2}(\sum_{j=0}^{2}(D+T)^{j})=D^{2}+\beta T,$$  
$$\pi_{3}(\sum_{j=0}^{3}(D+T)^{j})=D^{3}+\beta
DT+2(\beta-1)TD=D^{3}+(3\beta-2)TD+\beta(\partial T).$$
We can then write the formula at the third order :
$$\langle
S,(\nabla_{s}^{3}-(3(m+1)\lambda+2)r\vee\nabla_{s}-\lambda(m+1)(\nabla_{s}r))f\rangle$$
$$+C_{3,1}\langle
Div S,(\nabla_{s}^{2}-\lambda(m+1)r)f\rangle+C_{3,2}\langle(Div^{2}+\gamma_{5}(m+1)i(r))S,
\nabla_{s}f\rangle$$
$$+C_{3,3}\langle(Div^{3}+(3\gamma_{5}(m+1)-2)i(r)Div+\gamma_{5}(m+1)i(\nabla_{s}r))S,f\rangle.$$
At the second order, the formula is simply :
$$\langle
S,(\nabla_{s}^{2}-\lambda(m+1)r)f\rangle+C_{2,1}\langle
Div S,\nabla_{s}f\rangle+C_{2,2}\langle(Div^{2}+\gamma_{3}(m+1)i(r))S,f\rangle.$$
\section{Acknowledgements}
It is a pleasure to thank P. Mathonet for numerous fruitful discussions and for
his interest in our work. We thank the Belgian FRIA for his Research Fellowship.
\bibliographystyle{plain} \bibliography{formules}
\end{document}